\documentclass{amsart}

\usepackage[latin1]{inputenc}

\usepackage{amssymb}
\usepackage{amsmath}

\def\kro#1#2{\left(\frac{#1}{#2}\right)}
\def\<#1>{{\left\langle{#1}\right\rangle}}
\def\set#1{{\left\{{\def\st{\;:\;}#1}\right\}}}
\def\abs#1{{\left|{#1}\right|}}

\def\ZZ{\mathbb{Z}}

\begin{document}

\title[Computing central values of twisted $L$-series]
      {Computing central values of twisted $L$-series\\
       The case of composite levels}
\author{Ariel Pacetti}
\address{%
  Departamento de Matemática\\
  Universidad de Buenos Aires\\
  Pabellón I, Ciudad Universitaria. C.P:1428\\
  Buenos Aires, Argentina\\
}
\email{apacetti@dm.uba.ar}
\thanks{The first author was supported by a CONICET grant}
\author{Gonzalo Tornaría}
\address{%
  Centro de Matemática\\
  Facultad de Ciencias\\
  Iguá 4225 esq. Mataojo\\
  Montevideo, Uruguay\\
}
\email{tornaria@math.utexas.edu}
\maketitle

\section{Introduction}
Let $f \in S_2(N)$ be a newform of weight two and level $N$.
If $f(z) = \sum_{m=1}^\infty a(m)\,q^m$ where $q = e^{2 \pi iz}$,
and $D$ is a fundamental discriminant, we define the twisted
L-function
\[
   L(f,D,s) = \sum_{m=1}^\infty \frac{a(m)}{m^s} \kro{D}{m} \enspace.
\]
We will assume that the twisted $L$-series are primitive (i.e.
the corresponding twisted modular forms are newforms). There is no
loss of generality in making this assumption: if this is not the case,
then $f$ would be a quadratic twist of a newform of smaller level,
which we can choose instead.

The question of efficiently computing the family of central values $L(f,D,1)$,
for fundamental discriminants $D$, has been considered by several authors.
By Waldspurger's formula \cite{Waldspurger} these values are related to the Fourier
coefficients of certain modular forms of weight $3/2$.

In \cite{Gross}, Gross gives a
method to construct, for the case of prime level $p$, and provided $L(f,1) \neq 0$,
a weight $3/2$ modular form of level $4p$, and gives an explicit
version of Waldspurger's formula for the imaginary quadratic twists.
In \cite{BSP90} Böcherer and
Schulze-Pillot extend Gross's method to the case of square free
level, but their method works only for a fraction of imaginary
quadratic twists (determined by quadratic residue conditions).
Later in \cite{Pacetti-Tornaria} the case of level $p^2$ ($p$ a prime)
is considered, and this is used in \cite{Pacetti-Tornaria-proc},
provided $p\equiv 3\pmod{4}$, to compute
central values for \emph{real} quadratic twists.

In \cite{MRVT} the non-vanishing condition is removed, and two modular
forms of weight $3/2$ (one giving the imaginary quadratic twists and
another one giving the real quadratic twists) are constructed, in the
case of prime level.

The aim of this paper is to show how some of these ideas can be
combined to handle the case of composite levels.
In the case of odd squarefree level $N$, for instance, this method constructs $2^t$
modular forms, where $t$ is the number of prime factors of $N$, whose
coefficients give the central values of all the quadratic twists.
We will focus on examples for levels $N=27$, $N=15$, and $N=75$, which
exhibit most aspects of our methods.

\section{The curve $27A$}

Let $f$ be the modular form of level 27, corresponding to
the elliptic curve $X_0(27)$, of minimal equation
\[
    y^2+y=x^3-7 \enspace.
\]
The eigenvalue of $f$ for the Atkin-Lehner involution $W_{27}$ is $-1$, 
and the sign of the functional equation for $L(f,s)$ is $+1$.

Let $B=(-1,-3)$ be the quaternion algebra ramified at $3$ and
$\infty$, and consider the order $R = \<1,3i,\frac{1+3j}{2},\frac{i+k}{2}>$,
a Pizer order of reduced discriminant $27$. The class number of left
$R$-ideals for such order is $2$, and representatives for left
$R$-ideals are $\set{R,I}$ where $I =
\<4,12i,\frac{7+6i+3j}{2},\frac{6+13i+k}{2}>$.
The eigenvector for the Brandt matrices which 
corresponds to $f$ is $(1,-1)$, with height $3$.

The ternary quadratic forms associated to their right orders are
\begin{align*}
  Q_1(x,y,z) &= 4x^2+27y^2+28z^2-4xz \enspace, \\
\intertext{and}
 Q_2(x,y,z) &= 7x^2+16y^2+31z^2+16yz+2xz+4xy  \enspace,
\end{align*}
respectively.

Note that, since the twist of
$f$ by the quadratic character of conductor $3$ is $f$ itself,
we have
\[
  L(f,-3D,s) = L(f, D, s) \enspace,
\]
for $-3D$ a fundamental discriminant.
We will thus assume that $3\nmid D$.

\subsection{Imaginary quadratic twists}

Let $D<0$ be a fundamental discriminant. If $\kro{D}{3} = +1$, the
sign of the functional equation for $L(f,D,s)$ is $-1$, so its central
value vanishes trivially.  Hence we can
restrict to the case where $\kro{D}{3} = -1$.
In this case we can follow Gross's method, using classical theta
series
\[
   \Theta(Q_i) := \frac{1}{2} \sum_{(x,y,z)\in\ZZ^3} q^{Q_i(x,y,z)}\enspace;
\]
we obtain a weight $3/2$ modular form of level $4\cdot 27$,
namely
\[
   g = \Theta(Q_1) - \Theta(Q_2)
     = q^4 - q^7 - q^{19} + q^{28} - 2q^{40} + 2q^{43} + \cdots
     \enspace.
\]
\begin{table}
\begin{tabular}{||r|rr||r|rr||r|rr||}
\hline
$D$ & $c({D})$ & $L(f,D,1)$& $D$ & $c({D})$ &
$L(f,D,1)$& $D$ & $c({D})$ & $L(f,D,1)$\\
\hline
   -4 & 1 & 1.529954 & -67 & -1 & 0.373827 & -139 & 3 & 2.335842 \\
   -7 & -1 & 1.156537 & -79 & 1 & 0.344267 & -148 & 1 & 0.251523 \\
   -19 & -1 & 0.701991 & -88 & -2 & 1.304749 & -151 & -1 & 0.249012 \\
   -31 & 0 & 0.000000 & -91 & 1 & 0.320766 & -163 & -1 & 0.239670 \\
   -40 & -2 & 1.935256 & -103 & 1 & 0.301502 & -184 & 2 & 0.902318 \\
   -43 & 2 & 1.866526 & -115 & -2 & 1.141352 & -187 & -2 & 0.895051 \\
   -52 & 1 & 0.424333 & -127 & -2 & 1.086092 & -199 & -3 & 1.952200 \\
   -55 & 2 & 1.650392 & -136 & 2 & 1.049540 & & &  \\
\hline
\end{tabular}
\caption{Coefficients of $g$ and 
imaginary quadratic twists of $27A$ \label{table:27Ai}}
\end{table}
Table~\ref{table:27Ai} shows the values of the Fourier coefficients
$c({D})$ of $g$ and of $L(f,D,1)$, where $-200<D<0$ is a
fundamental discriminant such that $\kro{D}{3} = -1$.
The Gross type formula
\[
   L(f,D,1) = k\,\frac{\abs{c(D)}^2}{\sqrt{\abs{D}}}
   \enspace,\quad D<0\enspace,
\]
is satisfied, where $c(D)$ is the $\abs{D}$-th Fourier coefficient of $g$, and
\[
   k = \frac{1}{3}\cdot\frac{(f,f)}{L(f,1)} = 2 L(f,-4,1)\approx 3.059908074114385749826388345
   \enspace.
\]

\subsection{Real quadratic twists}
Let $D>0$ be a fundamental discriminant. In this case, if $\kro{D}{3}
= -1$ the sign of the functional equation for $L(f,D,s)$ will be $-1$,
and its central value will vanish trivially. For $\kro{D}{3} = +1$, we
will employ a method similar to the one used in \cite{MRVT} for prime
levels. We need to choose an auxiliary prime $l\equiv 3\pmod{4}$ such
that $\kro{-l}{3}=-1$ and such that $L(f,-l,1)\neq 0$, for example
$l=7$. Following \cite{MRVT} we define generalized theta series
\[
   \Theta_{-7}(Q_i) := \frac{1}{2} \sum_{(x,y,z)\in\ZZ^3}
   \omega^{(i)}_7(x,y,z)\;\omega^{(i)}_3(x,y,z)\,q^{Q_i(x,y,z)/7}\enspace,
\]
where $\omega_7$ and $\omega_3$ are the two kinds of weight function
introduced in \S2.2 and \S2.3 of \cite{MRVT}, respectively.
The superscript in $\omega^{(i)}_3$ and $\omega^{(i)}_7$ indicates
that we are writing the weight functions in the basis corresponding to
the quadratic form $Q_i$.

The weight function of the first kind can be computed as
\[
   \omega^{(1)}_7(x,y,z) = \begin{cases}
       0 & \text{if $7\nmid Q_1(x,y,z)$,} \\
       \kro{x}{7} & \text{if $7\nmid x$,} \\
       \kro{5z}{7} & \text{otherwise;}
   \end{cases}
\]
and
\[
   \omega^{(2)}_7(x,y,z) = \begin{cases}
       0 & \text{if $7\nmid Q_2(x,y,z)$,} \\
       \kro{3y+5z}{7} & \text{if $7\nmid 3y+5z$,} \\
       \kro{6x}{7} & \text{otherwise.}
   \end{cases}
\]

The weight function of the second kind can be computed as
\[
\omega^{(1)}_3(x,y,z) = \kro{x+z}{3},
\quad\text{and}\quad
\omega^{(2)}_3(x,y,z) = \kro{2x+y+2z}{3}.
\]

The generalized theta series will be
\[
  \Theta_{-7}(Q_1) = 
  -2q^{4}+2q^{13}+4q^{16}-4q^{25}+2q^{28}-2q^{37}-4q^{40}+\cdots\enspace,
\]
and
\[
  \Theta_{-7}(Q_2) =
  q-q^{4}-q^{13}+2q^{16}-3q^{25}-q^{28}+q^{37}+2q^{40}+\cdots\enspace.
\]
Note that
$\Theta_{-7}(Q_1)+2\Theta_{-7}(Q_2)
=2q-4q^{4}+8q^{16}-10q^{25}+\cdots$,
corresponding to the Eisenstein eigenvector for the Brandt matrices,
has nonzero Fourier coefficients only at square indices.
Since $\Theta_{-7}(Q_1)+2\Theta_{-7}(Q_2)\equiv
\Theta_{-7}(Q_1)-\Theta_{-7}(Q_2)\pmod{3}$,
this explains the fact
that the coefficients in Table~\ref{table:27Ar}, with the exception of
$c_{-7}(1)$, are all divisible by $3$.

Thus we obtain a modular form of weight $3/2$, namely
\[
   g_{-7} = \Theta_{-7}(Q_1) - \Theta_{-7}(Q_2)
      =
      q+q^{4}-3q^{13}-2q^{16}+q^{25}-3q^{28}+3q^{37}+6q^{40}+\cdots\enspace,
\]
and the formula is now
\[
   L(f,D,1) = k_{-7}\,\frac{\abs{c_{-7}(D)}^2}{\sqrt{\abs{D}}}
   \enspace,\quad D>0\enspace,
\]
where $c_{-7}(D)$ is the ${D}$-th Fourier coefficient of $g_{-7}$, and
\[
   k_{-7} = \frac{1}{3}\cdot\frac{(f,f)}{L(f,-7,1)\sqrt{7}}
          = L(f,1)\approx 0.5888795834284833191045631668
   \enspace .
\]
\begin{table}
\begin{tabular}{||r|rr||r|rr||r|rr||}
\hline
$D$ & $c_{-7}({D})$ & $L(f,D,1)$& $D$ & $c_{-7}({D})$ &
$L(f,D,1)$& $D$ & $c_{-7}({D})$ & $L(f,D,1)$\\
\hline
   1 & 1 & 0.588880 & 76 & -3 & 0.607942 & 136 & 6 & 1.817856 \\
   13 & -3 & 1.469932 & 85 & 0 & 0.000000 & 145 & 6 & 1.760536 \\
   28 & -3 & 1.001590 & 88 & -6 & 2.259892 & 157 & 6 & 1.691917 \\
   37 & 3 & 0.871301 & 97 & -3 & 0.538125 & 172 & 0 & 0.000000 \\
   40 & 6 & 3.351961 & 109 & 0 & 0.000000 & 181 & -9 & 3.545457 \\
   61 & 3 & 0.678585 & 124 & 6 & 1.903786 & 184 & -6 & 1.562860 \\
   73 & -3 & 0.620308 & 133 & -3 & 0.459561 & 193 & 3 & 0.381496 \\
\hline
\end{tabular}
\caption{Coefficients of $g_{-7}$ and 
real quadratic twists of $27A$ \label{table:27Ar}}
\end{table}
Table~\ref{table:27Ar} shows the values of the Fourier coefficients
$c_{-7}(D)$ of $g_{-7}$ and of $L(f,D,1)$, where $0<D<200$ is a
fundamental discriminant such that $\kro{D}{3} = 1$.

\section{The curve $15A$}

Let $f$ be the modular form of level $15$, corresponding to the
elliptic curve $X_0(15)$, of minimal equation
\[
   y^2 + xy + y = x^3 + x^2 - 10x - 10 \enspace.
\]
The eigenvalues of $f$ for the Atkin-Lehner involutions $W_3$ and
$W_5$ are $+1$ and $-1$, and the sign of the functional equation for
$L(f,s)$ is $+1$.

The method of Gross, as extended by Böcherer and Schulze-Pillot to the
case of squarefree levels, requires that the ramification of the
quaternion algebra agrees with the Atkin-Lehner eigenvalues. In this
case, it would be necessary to work with the quaternion algebra ramified at
$5$ and $\infty$. To exhibit the generality of our method, we will
work with the quaternion algebra ramified at $3$ and $\infty$ instead.

Let $B=(-1,-3)$ be such a quaternion algebra;
an Eichler order of level $15$ (index $5$ in a maximal order)
is given by $R = \<1,i,\frac{1+5j}{2},\frac{1+i+3j+k}{2}>$.
The number of clases of
left $R$-ideals is $2$, and a set of representatives of the classes is given by $\set{R,I}$
where $I = \<2,2i,\frac{3+2i+5j}{2},\frac{3+i+3j+k}{2}>$. The
eigenvector for the Brandt matrices corresponding to $f$ is $(1,-1)$,
with height $4$, and the ternary quadratic forms associated to $R$ and
$I$ are
\[
   Q_1(x,y,z) = Q_2(x,y,z) = 4x^2+15y^2+16z^2-4xz \enspace.
\]

\subsection{Imaginary quadratic twists}

Let $D<0$ be a fundamental discriminant. We say that $D$ is of type
$(s_1,s_2)$ if $\kro{D}{3}=s_1$ and $\kro{D}{5}=s_2$.
We need the sign of the functional equation for $L(f,D,s)$ to be $+1$,
so that its central value does not vanish trivially. For this to hold
we need $D$ to be of type $(-,+)$, $(+,-)$, $(+,0)$, $(0,-)$, or
$(0,0)$.

Note that the linear combination of classical theta series
$\Theta(Q_1)-\Theta(Q_2)$ is trivially zero, since $Q_1=Q_2$; this
reflects the fact that the ramification does not match the
Atkin-Lehner eigenvalues. Instead we set
\[
   \Theta_{1}(Q_i) := \frac{1}{4} \sum_{(x,y,z)\in\ZZ^3}
   \omega^{(i)}_3(x,y,z)\;\omega_5^{(i)}(x,y,z)\,q^{Q_i(x,y,z)}\enspace,
\]
where $\omega_3$ and $\omega_5$ are weight functions of the second
kind as in \cite[\S2.3]{MRVT}.
We have $\Theta_1(Q_1)=-\Theta_1(Q_2)$, and hence
we obtain a modular form of weight $3/2$ and level
$4\cdot 15^2$, namely
\[
   g_1 = 2\,\Theta_1(Q_1)
       = q^{4}+q^{16}+2q^{19}+2q^{31}+q^{64}+\cdots\enspace.
\]
The corresponding formula is
\[
   L(f,D,1) = k_{1}\,\frac{\abs{c_{1}(D)}^2}{\sqrt{\abs{D}}}
   \enspace,\quad\text{$D<0$ of type $(-,+)$,}
\]
where $c_1(D)$ is the $\abs{D}$-th Fourier coefficient of $g_1$, and
\[
   k_1 = \frac{1}{4}\cdot\frac{(f,f)}{L(f,1)}
          = 2L(f,-4,1)\approx 3.192484444263567020297938143
   \enspace ,
\]
\begin{table}
\begin{tabular}{||r|rr||r|rr||r|rr||}
\hline
$D$ & $c_{1}({D})$ & $L(f,D,1)$& $D$ & $c_{1}({D})$ &
$L(f,D,1)$& $D$ & $c_{1}({D})$ & $L(f,D,1)$\\
\hline
   -4 & 1 & 1.596242 & -91 & -4 & 5.354613 & -184 & -4 & 3.765649 \\
   -19 & 2 & 2.929625 & -136 & -4 & 4.380053 & -199 & -2 & 0.905237 \\
   -31 & 2 & 2.293549 & -139 & -2 & 1.083132 & & &  \\
   -79 & -2 & 1.436730 & -151 & 2 & 1.039203 & & &  \\
\hline
\hline
$D$ & $c_{17}({D})$ & $L(f,D,1)$& $D$ & $c_{17}({D})$ &
$L(f,D,1)$& $D$ & $c_{17}({D})$ & $L(f,D,1)$\\
\hline
   -3 & 2 & 0.921591 & -83 & 4 & 0.350421 & -152 & 8 & 1.035779 \\
   -8 & -4 & 1.128714 & -87 & 4 & 0.684541 & -155 & 8 & 2.051412 \\
   -15 & -2 & 0.824296 & -95 & 0 & 0.000000 & -167 & 4 & 0.247042 \\
   -20 & 4 & 1.427722 & -107 & 4 & 0.308629 & -168 & 8 & 1.970444 \\
   -23 & 4 & 0.665679 & -120 & -4 & 1.165730 & -183 & 0 & 0.000000 \\
   -35 & -4 & 1.079257 & -123 & -4 & 0.575713 & -195 & -4 & 0.914474
\\
   -47 & -4 & 0.465672 & -132 & -8 & 2.222961 & & &  \\
   -68 & 0 & 0.000000 & -143 & -8 & 1.067876 & & &  \\
\hline
\end{tabular}
\caption{Coefficients of $g_1$ and $g_{17}$, and imaginary twists of
$15A$\label{table:15Ai}}
\end{table}
c.f. Table \ref{table:15Ai} (top).

To obtain the other $4$ types of negative $D$, we need to choose an
auxiliary prime $l\equiv 1\pmod{4}$ such that
$\kro{l}{3}=\kro{l}{5}=-1$, and such that $L(f,l,1)\neq 0$,
e.g.  $l=17$. We then define the generalized theta series
\[
   \Theta_{17}(Q_i) := \frac{1}{4} \sum_{(x,y,z)\in\ZZ^3}
   \omega^{(i)}_{17}(x,y,z)\,q^{Q_i(x,y,z)/17}\enspace,
\]
where $\omega_{17}$ is the weight function of the first kind defined
in \cite[\S2.2]{MRVT}. Now
\[
   g_{17} = 2\,\Theta_{17}(Q_1) =
            2q^{3}-4q^{8}-2q^{15}+4q^{20}+4q^{23}+\cdots
\]
is a weight $3/2$ modular form of level $4\cdot 15$. As expected by
the multiplicity one theorem of Kohnen \cite{Kohnen}, this
form turns out to be the same as the one constructed by Böcherer and
Schulze-Pillot.
The formula in this case is
\[
   L(f,D,1) = \star\, k_{17}\,\frac{\abs{c_{17}(D)}^2}{\sqrt{\abs{D}}}
   \enspace,\quad\text{$D<0$ of type $(+,-)$, $(+,0)$, $(0,-)$, or
   $(0,0)$\enspace,}
\]
and $\star=1$, $2$, $2$, or $4$ respectively;
where $c_{17}(D)$ is the $\abs{D}$-th Fourier coefficient of $g_{17}$,
and
\[
   k_{17} = \frac{1}{4}\cdot\frac{(f,f)}{L(f,17,1)\sqrt{17}}
          \approx 0.1995302777664729387686211340
   \enspace ,
\]
c.f. Table \ref{table:15Ai} (bottom).

\subsection{Real quadratic twists}

Let $D>0$ be a fundamental discriminant. In order for the sign of the
functional equation of $L(f,D,s)$ to be $+1$, we need $D$ to be of
type $(+,+)$, $(0,+)$, $(-,-)$, or $(-,0)$.

For the first two types we need an auxiliary prime $l\equiv 3\pmod{4}$
such that $\kro{-l}{3}=-1$ and $\kro{-l}{5}=+1$, and such that
$L(f,-l,1)\neq 0$, e.g. $l=19$. Again
\[
   \Theta_{-19}(Q_i) := \frac{1}{4} \sum_{(x,y,z)\in\ZZ^3}
   \omega^{(i)}_{19}(x,y,z)\;\omega^{(i)}_{5}(x,y,z)\,q^{Q_i(x,y,z)/19}\enspace,
\]
with $\omega_{19}$ of the first kind and $\omega_{5}$ of the second
kind. The modular form
\[
   g_{-19} = 2\,\Theta_{-19}(Q_1) =
            2q-4q^{4}+2q^{9}-8q^{21}+8q^{24}+\cdots
\]
has level $4\cdot 15\cdot 5$, and the formula is
\[
   L(f,D,1) = \star\,k_{-19}\,\frac{\abs{c_{-19}(D)}^2}{\sqrt{\abs{D}}}
   \enspace,\quad\text{$D>0$ of type $(+,+)$ or $(0,+)$\enspace,}
\]
$\star=1$ or $2$ respectively;
$c_{-19}(D)$ is the ${D}$-th Fourier coefficient of $g_{-19}$, and
\[
   k_{-19} = \frac{1}{4}\cdot\frac{(f,f)}{L(f,-19,1)\sqrt{19}}
          = \frac{1}{4} L(f,1) \approx 0.08753769014578762644876130241
   \enspace .
\]
\begin{table}
\begin{tabular}{||r|rr||r|rr||r|rr||}
\hline
$D$ & $c_{-19}({D})$ & $L(f,D,1)$& $D$ & $c_{-19}({D})$ &
$L(f,D,1)$& $D$ & $c_{-19}({D})$ & $L(f,D,1)$\\
\hline
   1 & 2 & 0.350151 & 76 & -16 & 2.570563 & 141 & -8 & 0.943616 \\
   21 & -8 & 2.445093 & 109 & 16 & 2.146455 & 156 & 16 & 3.588416 \\
   24 & 8 & 2.287175 & 124 & 16 & 2.012446 & 181 & 0 & 0.000000 \\
   61 & 16 & 2.869261 & 129 & -8 & 0.986530 & 184 & -16 & 1.652061 \\
   69 & -8 & 1.348902 & 136 & 0 & 0.000000 & & &  \\
\hline
\end{tabular}
\begin{tabular}{||r|rr||r|rr||r|rr||}
\hline
$D$ & $c_{-23}({D})$ & $L(f,D,1)$& $D$ & $c_{-23}({D})$ &
$L(f,D,1)$& $D$ & $c_{-23}({D})$ & $L(f,D,1)$\\
\hline
   5 & 2 & 1.252737 & 77 & -8 & 2.553816 & 152 & 0 & 0.000000 \\
   8 & -4 & 1.980752 & 92 & 8 & 2.336367 & 173 & -12 & 3.833492 \\
   17 & 4 & 1.358785 & 113 & -4 & 0.527031 & 185 & 4 & 0.823795 \\
   53 & 4 & 0.769550 & 137 & -4 & 0.478646 & 188 & -8 & 1.634392 \\
   65 & -4 & 1.389787 & 140 & 8 & 3.787922 & 197 & 12 & 3.592398 \\
\hline
\end{tabular}
\caption{Coefficients of $g_{-19}$ and $g_{-23}$, and real twists of
$15A$\label{table:15Ar}}
\end{table}
Table \ref{table:15Ar} (top) shows the values of the coefficients
$c_{-19}(D)$ and the central values $L(f,D,1)$ for $0<D<200$ a
fundamental discriminant of type $(+,+)$ or $(0,+)$.

For the remaining two types we need an auxiliary prime $l\equiv 3\pmod{4}$
such that $\kro{-l}{3}=+1$ and $\kro{-l}{5}=-1$, and such that
$L(f,-l,1)\neq 0$, e.g. $l=23$. As before we define
\[
   \Theta_{-23}(Q_i) := \frac{1}{4} \sum_{(x,y,z)\in\ZZ^3}
   \omega^{(i)}_{23}(x,y,z)\;\omega^{(i)}_{3}(x,y,z)\,q^{Q_i(x,y,z)/23}\enspace,
\]
with $\omega_{23}$ of the first kind and $\omega_{3}$ of the second
kind. The modular form
\[
   g_{-23} = 2\,\Theta_{-23}(Q_1) =
   2q^{5}-4q^{8}+4q^{17}-4q^{32}+4q^{53}+\cdots
\]
has level $4\cdot 15\cdot 3$, and the formula is
\[
   L(f,D,1) = \star\,k_{-23}\,\frac{\abs{c_{-23}(D)}^2}{\sqrt{\abs{D}}}
   \enspace,\quad\text{$D>0$ of type $(-,-)$ or $(-,0)$\enspace,}
\]
$\star=1$ or $2$ respectively;
$c_{-23}(D)$ is the ${D}$-th Fourier coefficient of $g_{-23}$ and
\[
   k_{-23} = \frac{1}{4}\cdot\frac{(f,f)}{L(f,-23,1)\sqrt{23}}
          \approx 0.3501507605831505057950452092
   \enspace .
\]
Table \ref{table:15Ar} (bottom) shows the values of the coefficients
$c_{-19}(D)$ and the central values $L(f,D,1)$ for $0<D<200$ a
fundamental discriminant of type $(-,-)$ or $(-,0)$.

\section{The curve $75A$}

Let $f$ be the modular form of level 75 corresponding to
the elliptic curve of minimal equation
\[
    y^2+y=x^3-x^2-8x-7 \enspace.
\]
The eigenvalue of $f$ for the Atkin-Lehner involution $W_{3}$ is $+1$,
for $W_{25}$ is $-1$, and the sign of the functional equation for
$L(f,s)$ is $+1$.

Let $B=(-1,-3)$ be the quaternion algebra ramified at $3$ and
$\infty$, and consider the order $R = \<1,i,\frac{1+5j}{2},\frac{i+5k}{2}>$,
an Eichler order of level $75$ (index $25$ in a maximal order).
The class number of left $R$-ideals is $6$,
and the eigenvector for the Brandt matrices which 
corresponds to $f$ is $(1, -1, 1, -1, 0, 0)$, with height $6$.

The ternary quadratic forms associated to the right orders of the
choosen ideal class representatives are
\begin{align*}
  Q_1(x,y,z) = Q_2(x,y,z) & = 4x^2+75y^2+76z^2-4xz \enspace, \\
  Q_3(x,y,z) = Q_4(x,y,z) & = 16x^2+19y^2+79z^2+4xy+16xz+2yz \enspace, \\
\intertext{and}
  Q_5(x,y,z) = Q_6(x,y,z) & = 24x^2+31y^2+39z^2+24xy+12xz+6yz \enspace,
\end{align*}
respectively.

We will assume that $5\nmid D$. Indeed, the twist of $f$ by the
quadratic character of conductor $5$ is another modular form $f'$ of
level $75$, thus we have
\[
   L(f,5D,1) = L(f',D,1),
\]
for $5D$ a fundamental discriminant. By applying the same procedure to
the modular form $f'$ we can compute the central values for these
twists. So, we actually need $8$ different modular forms of weight
$3/2$ to compute all the twisted central values.

\subsection{Imaginary quadratic twists}

Let $D<0$ be a fundamental discriminant. If the sign of the
functional equation for $L(f,D,s)$ is $+1$, the type of $D$ has to be
either $(-,+)$ or $(-,-)$.

For the first case we look at the generalized theta
series
\[
   \Theta_1(Q_i) := \frac{1}{4} \sum_{(x,y,z)\in\ZZ^3}
   \omega^{(i)}_3(x,y,z)\;\omega^{(i)}_5(x,y,z)\,q^{Q_i(x,y,z)}\enspace;
\]
we obtain the modular form
\[
   g_1 = 2\Theta_1(Q_1) - 2\Theta_1(Q_3)
     = q^{4}-2q^{16}-q^{19}-q^{31}-2q^{64}+3q^{76}+4q^{79}-q^{91}+\cdots
     \enspace.
\]
The formula
\[
   L(f,D,1) = k_1\,\frac{\abs{c_1(D)}^2}{\sqrt{\abs{D}}}
   \enspace, \quad\text{$D<0$ of type $(-,+)$\enspace,}
\]
is satisfied (c.f. Table \ref{table:75Ai}, top),
where $c_1(D)$ is the $\abs{D}$-th Fourier coefficient of $g_1$ and
\[
   k_1 = \frac{1}{6}\cdot\frac{(f,f)}{L(f,1)}
   = 2\,L(f,-4,1)\approx 4.669532748718719327951206761
   \enspace.
\]

In the second case we need to choose an auxiliary prime $l \equiv 1
\pmod{4}$ such that $\kro{l}{3} = +1$, $\kro{l}{5} = -1$ and $L(f,l,1)
\neq 0$, for example $l=13$, and define
\[
   \Theta_{13}(Q_i) := \frac{1}{4} \sum_{(x,y,z)\in\ZZ^3}
   \omega^{(i)}_{13}(x,y,z)\;\omega^{(i)}_3(x,y,z)\;\omega^{(i)}_5(x,y,z)
      \,q^{Q_i(x,y,z)/13}\enspace.
\]
We obtain the modular form
\[
   g_{13} = 2\Theta_{13}(Q_1) - 2\Theta_{13}(Q_3)
     = 3q^{7}+3q^{28}+3q^{43}+3q^{52}-3q^{67}-6q^{88}+\cdots
     \enspace,
\]
and the formula
\[
   L(f,D,1) = k_{13}\,\frac{\abs{c_{13}(D)}^2}{\sqrt{\abs{D}}}
   \enspace,\quad\text{$D<0$ of type $(+,+)$\enspace,}
\]
is satisfied (c.f. Table \ref{table:75Ai}, bottom),
where $c_{13}(D)$ is the $\abs{D}$-th Fourier coefficient of $g_{13}$ and
\[
   k_{13} = \frac{1}{6}\cdot\frac{(f,f)}{L(f,13,1)\sqrt{13}}
   \approx 1.556510916239573109317068920
   \enspace.
\]

\begin{table}
\begin{tabular}{||r|rr||r|rr||r|rr||}
\hline
$D$ & $c_{1}({D})$ & $L(f,D,1)$& $D$ & $c_{1}({D})$ &
$L(f,D,1)$& $D$ & $c_{1}({D})$ & $L(f,D,1)$\\
\hline
   -4 & 1 & 2.334766 & -91 & -1 & 0.489500 & -184 & 2 & 1.376970 \\
   -19 & -1 & 1.071264 & -136 & 2 & 1.601637 & -199 & -5 & 8.275360 \\
   -31 & -1 & 0.838673 & -139 & -2 & 1.584258 & & &  \\
   -79 & 4 & 8.405816 & -151 & 5 & 9.500030 & & &  \\
\hline
\hline
$D$ & $c_{13}({D})$ & $L(f,D,1)$& $D$ & $c_{13}({D})$ &
$L(f,D,1)$& $D$ & $c_{13}({D})$ & $L(f,D,1)$\\
\hline
   -7 & 3 & 5.294752 & -88 & -6 & 5.973286 & -163 & 3 & 1.097238 \\
   -43 & 3 & 2.136291 & -103 & -6 & 5.521233 & -187 & 0 & 0.000000 \\
   -52 & 3 & 1.942643 & -127 & -6 & 4.972248 & & &  \\
   -67 & -3 & 1.711423 & -148 & 0 & 0.000000 & & &  \\
\hline
\end{tabular}
\caption{Coefficients of $g_1$ and $g_{13}$, and imaginary twists of
$75A$\label{table:75Ai}}
\end{table}

\subsection{Real quadratic twists}
Let $D>0$ be a fundamental discriminant. The only possibilities so
that the sign of the functional equation for $L(f,D,s)$ is $+1$ are
the discriminants $D$ of types $(+,+), (0,+), (+,-)$, and $(0,-)$.

For the first two cases we can use the generalized theta series
\[
   \Theta_{-19}(Q_i) := \frac{1}{2} \sum_{(x,y,z)\in\ZZ^3}
   \omega^{(i)}_{19}(x,y,z)\;\omega^{(i)}_5(x,y,z)\, q^{Q_i(x,y,z)/19}\enspace.
\]
Thus we obtain a modular form of weight $3/2$, namely
\[
   g_{-19} = q+q^{4}+q^{9}-q^{21}-2q^{24}-q^{36}-4q^{49}
       -q^{61}+\cdots
   \enspace,
\]
and the formula is
\[
   L(f,D,1) = \star\, k_{-19}\,\frac{\abs{c_{-19}(D)}^2}{\sqrt{\abs{D}}}
   \enspace,\quad\text{$D>0$ of type $(+,+)$ or $(0,+)$\enspace,}
\]
$\star=1$ or $2$ respectively, $c_{-19}(D)$ the ${D}$-th Fourier
coefficient of $g_{-19}$, and
\[
   k_{-19} = \frac{1}{6}\cdot\frac{(f,f)}{L(f,-19,1)\sqrt{19}}
          = L(f,1)\approx 1.402539940216221119844494086
   \enspace ,
\]
c.f. Table \ref{table:75Ar} (top).

\begin{table}
\begin{tabular}{||r|rr||r|rr||r|rr||}
\hline
$D$ & $c_{-19}({D})$ & $L(f,D,1)$& $D$ & $c_{-19}({D})$ &
$L(f,D,1)$& $D$ & $c_{-19}({D})$ & $L(f,D,1)$\\
\hline
   1 & 1 & 1.402540 & 76 & 1 & 0.160882 & 141 & 2 & 0.944921 \\
   21 & -1 & 0.612119 & 109 & -1 & 0.134339 & 156 & -1 & 0.224586 \\
   24 & -2 & 2.290338 & 124 & 5 & 3.148795 & 181 & 3 & 0.938250 \\
   61 & -1 & 0.179577 & 129 & 5 & 6.174338 & 184 & -2 & 0.413586 \\
   69 & 2 & 1.350768 & 136 & -6 & 4.329605 & & &  \\
\hline
\hline
   12 & 3 & 2.429270 & 73 & 6 & 1.969859 & 168 & 6 & 2.596999 \\
   13 & 3 & 1.166984 & 88 & -6 & 1.794135 & 172 & 3 & 0.320828 \\
   28 & -3 & 0.795165 & 93 & -3 & 0.872620 & 177 & -6 & 2.530113 \\
   33 & -6 & 5.859621 & 97 & 9 & 3.844972 & 193 & 9 & 2.725840 \\
   37 & 0 & 0.000000 & 133 & -3 & 0.364847 & & &  \\
   57 & -3 & 1.114626 & 157 & 3 & 0.335805 & & &  \\
\hline
\end{tabular}
\caption{Coefficients of $g_{-19}$ and $g_{-7}$, and real twists of
$75A$\label{table:75Ar}}
\end{table}

In the other two cases we can use the generalized theta series
\[
   \Theta_{-7}(Q_i) := \frac{1}{2} \sum_{(x,y,z)\in\ZZ^3}
   \omega^{(i)}_{7}(x,y,z)\;\omega^{(i)}_5(x,y,z)\, q^{Q_i(x,y,z)/7}\enspace.
\]
We obtain a modular form of weight $3/2$
\[
   g_{-7} = 3q^{12}+3q^{13}-3q^{28}-6q^{33}+6q^{48}-9q^{52}-3q^{57}+6q^{73}+\cdots
     \enspace,
\]
satisfying the formula
\[
   L(f,D,1) = \star\,k_{-7}\,\frac{\abs{c_{-7}(D)}^2}{\sqrt{\abs{D}}}
   \enspace,\quad\text{$D>0$ of type $(+,-)$ or $(0,-)$\enspace,}
\]
$\star=1$ or $2$ respectively, $c_{-7}(D)$ the ${D}$-th Fourier
coefficient of $g_{-7}$, and
\[
   k_{-7} = \frac{1}{6}\cdot\frac{(f,f)}{L(f,-7,1)\sqrt{7}}
          \approx 0.4675133134054070399481646950
   \enspace .
\]
c.f. Table \ref{table:75Ar} (bottom).

\end{document}